\documentclass{elsarticle}

\usepackage{graphicx}      
\usepackage{geometry}
\usepackage{multicol}
\usepackage{enumerate}
\usepackage{amsmath}
\usepackage{xcolor}
\usepackage{amsfonts}
\usepackage{float}
\usepackage{MnSymbol}%
\usepackage{wasysym}%
\usepackage{algorithm}
\graphicspath{{Figures/}}

\newtheorem{remark}{Remark}

\DeclareMathOperator*{\argmin}{arg\,min}

\def\R{\mathbb{R}}
\def\L{\mathbb{L}}
\def\Id{{\rm{I}}}

\begin{document}

\title{Modifier-Adaptation for Real-Time Optimal Periodic Operation} 

\author[1]{Victor Mirasierra} 
\author[1]{Daniel Limon}

\date{}

\journal{ }

\address[1]{Departamento de Ingenier\'ia de Sistemas y Autom\'atica, Universidad de Sevilla, Escuela Superior de Ingenieros, CO 41092 Spain. (e-mail: \{vmirasierra,dlm\}@us.es).}

\begin{abstract}                
In this paper, we present the periodic modifier-adaptation formulation of the dynamic real time optimization. The proposed formulation uses gradient information to update the problem with affine modifiers so that, upon convergence, its solution matches the optimal steady periodic trajectory. Unlike other state of the art modifier-adaptation techniques, the proposed approach is able to converge not only to optimal steady states, but also to optimal periodic trajectories. The full control scheme to take the system from its current state to the optimal periodic trajectory is detailed. The convergence of the computed reference to the optimal periodic behaviour is shown by means of a periodic version of the quadruple tank benchmark.
\end{abstract}

\maketitle

\section{Introduction}

Economic optimization plays a major role in most industries, since it allows to optimize the performance of the real plant operation \cite{singh2012overview}. To achieve optimal performance, optimization problems leverage system data and models to compute the trajectory that minimize the economic cost. To arrange the system information into an optimization framework, multidisciplinary teams are often involved. They must have deep knowledge about the real systems and be able to build detailed models that mirror their behaviour. The complexity and possible change over time of real systems (e.g. due to deterioration) make the identification task expensive and prone to errors, which leads to plant-model mismatch and ultimately may lead to a loss of performance in the controlled system.

In two-layer control schemes \cite{findeisen1980control}, the economic optimization is splitted in two main layers. The first one, called real time optimization (RTO), computes the optimal steady behaviour, while the second one, known as advanced control, calculates the inputs required to take the system to that reference. To transform the optimal reference computed by the RTO into a valid reference to the advance controller, often an intermediate layer known as the steady-state target optimization (SSTO) \cite{muske1997steady} is used.

One of the strengths of two-layer control schemes is their ability to use different models for the different layers, being the model from the RTO layer usually more complex and global, while the one from the advanced control layer is generally faster and able to quickly react to disturbances. This allows to keep a high performance from the detailed RTO, while keeping the control fast from the advanced control.

Standard formulations of the RTO deal with the optimization of the plant operated at equilibrium points. However, there exist many scenarios where the plant operates optimally with a periodic behaviour, such as HVAC systems, solar plants, water distribution networks, electric networks, among others. For these systems with periodic nature, a dynamic RTO is better suited because of its ability to calculate not only the optimal steady state, but also the optimal periodic trajectory. This constitutes a generalization of the standard RTO and usually comes at the expense of an increased complexity because of the larger number of variables. While dynamic RTO schemes may theoretically converge to the optimal steady operation, it is still sensitive to plant-model mismatch, thus making it susceptible to a performance decrease.

In order to cope with the issues derived from the plant-model mismatch, modifier-adaptation (MA) formulations of the RTO emerged and have been studied over the last decades \cite{marchetti2016modifier,mirasierra2020real,rodriguez2018optimizacion,vergara2020modifier} with promising results. They update the model-based RTO problem with affine modifiers that incorporate information of the real system. Upon convergence of the modifiers, the modified problem is able to calculate either the optimal steady operation of the real system or the optimal input profile of a batch process \cite{costello2011modifier} from an initially inaccurate model. Modifier-adaptation schemes have been mainly built upon the standard RTO to compute the optimal steady state of a system. In this work we present a periodic modifier adaptation scheme which is built upon a dynamic RTO and calculates upon convergence the optimal periodic trajectory of a real system. The proposed approach can be seen as a generalization of the MA scheme proposed in \cite{marchetti2009modifier} to include optimal periodic behaviour.

The structure of the paper is the following: In Section \ref{sec:problem:formulation} we introduce the problem under consideration, along with the two-layer control scheme. Then, in Section \ref{sec:pma} we analyse how to modify the dynamic RTO so that, upon convergence, its solution matches the first order necessary conditions of optimality of the optimal problem. Then, in Section \ref{sec:stto} we detail how to transform the optimal operation computed by the dynamic RTO into a valid reference for the advanced control layer. Section \ref{sec:mpc} shows a way to design the advanced layer to follow a reference. In Section \ref{sec:algorithm} we present the full algorithms to implement the two-layer control scheme with periodic modifier-adaptation. A simplified version of this algorithm is used in Section \ref{sec:example} on the quadruple tank benchmark example to test the performance of the proposed approach. Finally, Section \ref{sec:conclusions} discusses the conclusions.

\section{Problem formulation} \label{sec:problem:formulation}
Consider a system that is described by the following (unknown) discrete-time state-space representation:

\begin{equation} \label{eq:real:sys}
x_{k+1} = f_{p,k}(x_k,u_k),
\end{equation}
where $x_k \in \R^{n_x}$ and $u_k \in \R^{n_u}$ are respectively the states and inputs of the system at time $k$, and $f_{p,k} : \R^{n_x \times n_u} \to \R^{n_x}$ represents the dynamics of the real system at time $k$. Each step in $k$ represents $t_T$ seconds.

Let system \eqref{eq:real:sys} be periodic with known period $T t_T$ seconds, i.e. $f_{p,k} = f_{p,k+T}$, and let $x_0$ be the initial state. At the first step of each period, given the sequence of $T$ next inputs ${\bf{u}}_T = \begin{bmatrix} u_0^T & u_1^T & \cdots & u_{T-1}^T \end{bmatrix}^T \in \R^{T n_u}$, then the $T$ following states of the system \eqref{eq:real:sys} are defined by the time-invariant function $F_p : \R^{n_x \times T n_u} \to \R^{T n_x}$ so that:

\begin{equation} \label{eq:SYS}
{\bf{x}}_T = \begin{bmatrix} x_{1}^T & x_{2}^T & \cdots & x_{T}^T \end{bmatrix}^T = F_{p}(x_0,{\bf{u}}_T).
\end{equation}

At any time $k$, the states and inputs of system \eqref{eq:real:sys} can be subject to (possibly nonlinear) constraints of the form:

\begin{equation} \label{eq:SYS:cons}
g_k(x_k,u_k) \leq 0,
\end{equation}
which are also periodic with period $T t_T$ seconds. Considering the periodic constraint $x_0 = x_T$, at the first step of each period the constraints \eqref{eq:SYS:cons} can also be expressed by its compact form:
$$G({\bf{x}}_T,{\bf{u}}_T) \leq 0,$$
where $G : \R^{T n_x \times T n_u} \to \R$.

The optimal economic control problem calculates the infinite sequence of inputs that, when applied to the system \eqref{eq:real:sys}, minimizes the economic cost given by the stage cost function $\phi_k(x_k,u_k)$ over time. Let the stage cost function $\phi_k$ be periodic with period $T t_T$ seconds and consider the periodic constraint $x_0 = x_T$, then at the first step of each period, the time-invariant cost function $\Phi$ represents the sum of stage cost functions $\phi_k$ over the $T$ future steps and is defined as:

\begin{equation}
\Phi({\bf{x}}_T,{\bf{u}}_T) = \sum_{i=0}^{T-1} \phi_{i}(x_i,u_i).
\end{equation}

Given the initial state of the system $x_0$, the optimal economic control problem can be formulated as follows:

\begin{equation} \label{eq:opt:control}
\begin{aligned}
\min_{{\bf{u}}_\infty} \quad &\sum_{k=0}^{\infty} \phi_k(x_k,u_k) \\
\text{s.t.} \quad & x_{k+1} = f_{p,k}(x_k,u_k), \quad \text{for all } k=0, 1, \ldots, \infty \\
& g_k(x_k,u_k)\leq 0, \quad \text{for all } k=0, 1, \ldots, \infty 
\end{aligned}
\end{equation}

In real applications, the previous fomulation is seldom implemented because of two main reasons: (i) the system dynamics (i.e. $f_{p,k}$) are usually unknown, and (ii) the infinite number of decision variables hinders the problem's readiness for implementation. 

\begin{figure}
\centering
\includegraphics[scale=0.5]{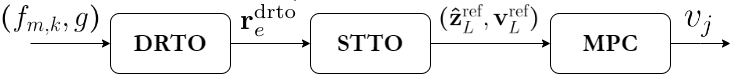}
\caption{Control diagram.}
\label{fig:stto}
\end{figure}

In practice, problem \eqref{eq:opt:control} is often tackled using a two-layer control scheme (Figure \ref{fig:stto}). In this scheme, the upper layer, also known as real time optimization (RTO), calculates the optimal operation of the system. Whereas the lower one, known as advanced control, computes the input sequence required to take the system from its current state to a given reference. These layers are separate and usually use different models and time horizons, so we consider also an intermediate layer called steady trajectory target optimization (STTO) which turns the optimal operation computed by the upper layer into a valid steady reference for the lower layer.

The basis of the two layer architecture is to split the main control problem into two smaller problems of different complexities and which are solved with different frequencies. On one hand, the RTO usually works with a complex and accurate model of the global plant. This model typically describes the fundamental and generally slow behaviour of the plant, which results in large time scales and low update frequency for the RTO. On the other hand, the advanced control generally uses a local dynamic model of system. It uses simple and fast models and its time scales are short. One of the benefits of this scheme is that the upper layer does not need to be recalculated with the same frequency as the lower one. This reduces computational costs, while keeping the control fast.

In this work we use dynamic real-time optimization (DRTO) as the upper layer. Unlike standard RTO, which aims to calculate the optimal steady setpoint ($x^s,u^s$), the objective of the DRTO is to compute the optimal periodic trajectory (${\bf{\hat{x}}}^{\text{drto}}_T, {\bf{u}}^{\text{drto}}_T$) with a predefined period of $T t_T$ seconds. The optimal periodic trajectory can be seen as a generalization of the optimal steady setpoint, since they lead to the same solution for $T=1$. Consequently, the DRTO can lead to better steady performance than the standard RTO, at the expense of it being a more intricate problem. In the case of periodic systems, it has been proven that the DRTO formulation is able to capture their optimal steady operation \cite{limon2014single}.

The DRTO uses a model of the real system $F_{m}$, instead of the real system dynamics $F_{p}$ described in \eqref{eq:SYS}:

\begin{equation} \label{eq:model:sys}
{\bf{\hat{x}}}_T = \begin{bmatrix} \hat{x}_{1}^T & \hat{x}_{2}^T & \cdots & \hat{x}_{T}^T \end{bmatrix}^T = F_{m}(x_0,{\bf{u}}_T),
\end{equation}
where $\hat{x}_{k}$ is the state predicted by the model at time $k$. Like \eqref{eq:SYS}, each step in $k$ equals $t_T$ seconds. Because of the complexity of real systems, models are usually unable to perfectly capture the real dynamics, leading to plant-model mismatch, i.e. $x_{k+1} \neq \hat{x}_{k+1}$.

One iteration of the DRTO is solved every $t_D = D (T t_T)$ seconds, with $D$ being a positive integer. Given the period $T$, they can be formulated as:

\begin{equation} \label{eq:drto}
\begin{aligned}
( x^{\text{drto}}_0, {\bf{\hat{x}}}^{\text{drto}}_T, &{\bf{u}}^{\text{drto}}_T ) = \\
\argmin_{x_0,{\bf{\hat{x}}}_T,{\bf{u}}_T} \quad &\sum_{i=0}^{T-1} \phi_i(\hat{x}_i,u_i) \\
\text{s.t.} \quad & {\bf{\hat{x}}}_T = F_{m}(x_0,{\bf{u}}_T) \\
& G({\bf{\hat{x}}}_T,{\bf{u}}_T) \leq 0 \\
& \hat{x}_T = x_0.
\end{aligned}
\end{equation}

The aforementioned formulation of the DRTO computes the optimal periodic operation for the available model of the system. However, due to plant-model mismatch, we know that this operation may not be optimal for the real system and might even lead to constraint violation. In the next section we present a reformulation of \eqref{eq:drto} which uses gradient based modifiers to update the base model so that, upon convergence, the solution of the modified DRTO matches the optimal periodic operation. Later, in Sections \ref{sec:mpc} and \ref{sec:stto}, the advanced control and the steady trajectory target optimization layers will be detailed.

\section{Periodic Modifier-Adaptation} \label{sec:pma}

Modifier-adaptation (MA) methodologies arose to correct the plant-model mismatch at the RTO level \cite{marchetti2016modifier,marchetti2009modifier}. They use measures and gradients from the system to build modifiers that update the RTO with affine terms. Upon convergence, MA schemes guarantee the satisfaction of the first order necessary conditions for optimality of the optimal problem. Traditionally, MA schemes have been built upon the standard RTO, which ultimately calculates the optimal steady setpoint. In this section we generalize state of the art approaches and show how to apply modifier-adaptation to the DRTO problem \eqref{eq:drto} and address the plant-model mismatch for optimal periodic trajectories. Zeroth and first order modifiers will be presented to update the dynamic model and ensure that, upon convergence, the optimal solution of the modified DRTO matches the optimal periodic trajectory of the system.

Let each iteration of the DRTO be labelled by index $l$. Then, given the modifiers $\lambda^x_l \in \R^{T n_x \times n_x}$, $\lambda^u_l \in \R^{T n_x \times T n_u}$ and $\epsilon_l \in \R^{T n_x}$, we introduce the periodic modifier-adaptation (P-MA) formulation of the DRTO at iteration $l$:

\begin{equation} \label{eq:drto:ma}
\begin{aligned}
( x^{\text{drto}}_0, {\bf{\hat{x}}}^{\text{drto}}_T, &{\bf{u}}^{\text{drto}}_T ) = \\
\argmin_{x_0,{\bf{\hat{x}}}_T,{\bf{u}}_T} \quad &\Phi({\bf{\hat{x}}}_T,{\bf{u}}_T) \\
\text{s.t.} \quad & {\bf{\hat{x}}}_T = F_m(x_0,{\bf{u}}_T) + \lambda^x_l x_0 + \lambda^u_l {\bf{u}}_T + \epsilon_l\\
& G({\bf{\hat{x}}}_T,{\bf{u}}_T)\leq 0 \\
& M {\bf{\hat{x}}}_T = x_0.
\end{aligned}
\end{equation}
where $M$ represents the constant matrix that ensures that the periodic constraint meets, i.e. $\hat{x}_T = x_0$.

After solving problem \eqref{eq:drto:ma}, the DRTO identifies a set of variables ${\bf{r}}_e^\text{drto} \in \R^{T n_r}$ that univocally defines the optimal economic trajectory 

\begin{equation} \label{eq:drto:re}
{\bf{r}}_e^{\text{drto}} = r_e( {\bf{\hat{x}}}^{\text{drto}}_T, {\bf{u}}^{\text{drto}}_T )
\end{equation}
and passes it to the STTO, which then transforms it into a valid reference for the MPC $({\bf{\hat{z}}}^{\text{ref}}_{N,j}, {\bf{v}}^{\text{ref}}_{N,j})$ (See Figure \ref{fig:stto}).

Now, we see how to calculate the modifiers $\lambda^x_l, \lambda^u_l$ and $\epsilon_l$ so that, upon convergence, the KKT conditions of problem \eqref{eq:drto:ma} converge to those of the optimal problem.

\subsection{KKT Matching} \label{sec:kkt}

In this section we show how to update the modifiers $\lambda^x_l, \lambda^u_l$ and $\epsilon_l$ so that the first order necessary conditions of optimality, also known as KKT conditions, of the P-MA DRTO \eqref{eq:drto:ma} match with those of the optimal problem. 

Given period $T$, the real optimal periodic trajectory $({\bf{x}}^{\text{opt}}_T,{\bf{u}}^{\text{opt}}_T)$ can be computed as the optimal solution to the following optimization problem:

\begin{equation} \label{eq:drto:opt}
\begin{aligned}
( x^{\text{opt}}_0, {\bf{x}}^{\text{opt}}_T, &{\bf{u}}^{\text{opt}}_T ) = \\
\argmin_{x_0,{\bf{x}}_T,{\bf{u}}_T} \quad &\Phi({\bf{x}}_T,{\bf{u}}_T) \\
\text{s.t.} \quad & {\bf{x}}_T = F_p(x_0,{\bf{u}}_T)\\
& G({\bf{x}}_T,{\bf{u}}_T)\leq 0 \\
& M {\bf{x}}_T = x_0.
\end{aligned}
\end{equation}

For the sake of simplicity and comparison, we define $\theta=\begin{bmatrix} x_0 \\ {\bf{u}}_T \end{bmatrix}$ and reformulate \eqref{eq:drto:opt} as:

\begin{subequations} \label{eq:kkt:opt}
\begin{align}
\begin{split}
\min_\theta \quad & \Phi^\theta(F_p^\theta(\theta),\theta) 
\end{split}\\
\begin{split}
\text{s.t.} \quad & G^\theta(F^\theta_p(\theta),\theta) \leq 0
\end{split} \\
\begin{split} \label{subeq:kkt:opt:periodic}
& M_1 F^\theta_p(\theta) + M_2 \theta = {\bf{0}} ,
\end{split}
\end{align}
\end{subequations}
where $\Phi^\theta, F_p^\theta$, $G^\theta$, $M1$ and $M2$ are functions derived from rewritting the ones in \eqref{eq:drto:opt} in terms of $\theta$, i.e. ${\bf{x}}_T = F_p^\theta(\theta)$ and similar, and \eqref{subeq:kkt:opt:periodic} corresponds to the periodic constraint.

We also define a modified version of the dynamic RTO \eqref{eq:drto:ma} at step $l$:
\begin{equation} \label{eq:kkt:ma}
\begin{aligned}
\min_{\theta} \quad & \Phi^\theta(F^\theta_m(\theta)+(\Lambda^\theta_l)^T \theta + \epsilon_l,\theta) \\
\text{s.t.} \quad & G^\theta(F^\theta_m(\theta)+(\Lambda^\theta_l)^T \theta + \epsilon_l,\theta) \leq 0 \\
 & M_1 (F_m^\theta(\theta)+(\Lambda^\theta_l)^T \theta + \epsilon_l) + M_2 \theta = {\bf{0}} \\
\end{aligned}
\end{equation}
where $\epsilon_l$ and $\Lambda^\theta_l = \begin{bmatrix} \lambda^x_l & \lambda^u_l	\end{bmatrix}^T$ refers to the zeroth and first order modifiers respectively at step $l$.

The Lagrangian function associated to the problem \eqref{eq:kkt:opt} is:

\begin{equation*}
\begin{aligned}
\L_p(\theta) = & \Phi^\theta(F_p^\theta(\theta),\theta) + \pi_1^T \left( G^\theta(F_p^\theta(\theta),\theta) \right) + \\
&\pi_2^T \left( M_1 F_p^\theta(\theta)+M_2 \theta \right),
\end{aligned}
\end{equation*}

and its gradient with respect to the decision variable $\theta$ is:

\begin{equation*}
\begin{aligned}
\frac{\partial \L_p}{\partial \theta} = & \frac{\partial \Phi^\theta}{\partial F^\theta_p} (\frac{\partial F^\theta_p}{\partial \theta}) + \frac{\partial \Phi^\theta}{\partial \theta} + \pi_1^T \Big[ \frac{\partial G^\theta}{\partial F^\theta_p} (\frac{\partial F^\theta_p}{\partial \theta}) + \frac{\partial G^\theta}{\partial \theta} \Big] + \\
& \pi_2^T \Big(M_1(\frac{\partial F^\theta_p}{\partial \theta})+M_2 \Big). \\
\end{aligned}
\end{equation*}

Analogously, the gradient of the Lagrangian function associated to problem \eqref{eq:kkt:ma} is the following:

\begin{equation*}
\begin{aligned}
\frac{\partial \L_m}{\partial \theta} = &\frac{\partial \Phi^\theta}{\partial F^\theta_m} (\frac{\partial F^\theta_m}{\partial \theta}+ \Lambda_\infty^\theta) + \frac{\partial \Phi^\theta}{\partial \theta} +\\
&\pi_1^T  \Big[ \frac{\partial G^\theta}{\partial F^\theta_m} (\frac{\partial F^\theta_m}{\partial \theta}+\Lambda_\infty^\theta) +  \frac{\partial G^\theta}{\partial \theta}  \Big] + \\
&\pi_2^T \Big(M_1(\frac{\partial F^\theta_m}{\partial \theta}+\Lambda_\infty^\theta)+ M_2 \Big).\\
\end{aligned}
\end{equation*}

Let $\theta^*$ be the (a priori unknown) optimal operation of the system, then the KKT conditions associated to problem \eqref{eq:kkt:opt} are:

\begin{subequations}
\begin{align}
\begin{split}
& \frac{\partial \L_p}{\partial \theta}(\theta^*) = 0 
\end{split} \\
\begin{split}
& G^\theta(F^\theta_p(\theta^*),\theta^*) \leq 0
\end{split} \\
\begin{split}
& M_1 F^\theta_p(\theta) + M_2 \theta = {\bf{0}}
\end{split} \\
\begin{split}
& \pi_1^*, \pi_2^* \geq 0
\end{split} \\
\begin{split}
& (G^\theta(F^\theta_p(\theta^*),\theta^*) ) \pi_1^* =0
\end{split} \\
\begin{split}
& (M_1 F^\theta_p(\theta) + M_2 \theta)_j \pi_{2,j}^* = 0, \; \; j=0,1, \ldots, n_x + T n_u.
\end{split}
\end{align}
\end{subequations}

Analogously, the KKT conditions associated to problem \eqref{eq:kkt:ma} are:

\begin{subequations}
\begin{align}
\begin{split}
& \frac{\partial \L_m}{\partial \theta}(\theta^*) = 0
\end{split} \\
\begin{split}
& G^\theta(F_m^\theta(\theta^*)+(\Lambda_l^\theta)^T \theta^* + \epsilon_l,\theta^*) \leq 0
\end{split} \\
\begin{split}
& M_1 (F_m^\theta(\theta^*)+(\Lambda_l^\theta)^T \theta^* + \epsilon_l) + M_2 \theta^* = {\bf{0}}
\end{split} \\
\begin{split}
& \pi_1^*, \pi_2^* \geq 0
\end{split} \\
\begin{split}
& (G^\theta(F_m^\theta(\theta^*)+(\Lambda_l^\theta)^T \theta^* + \epsilon_l,\theta^*) ) \pi_1^* =0 
\end{split} \\
\begin{split}
& \left( M_1 ( F_m^\theta(\theta^*)+(\Lambda_l^\theta)^T \theta^* + \epsilon_l ) + M_2 \theta^* \right)_j \pi_{2,j}^* = 0, \\ &j=0,1, \ldots, n_x + T n_u,
\end{split}
\end{align}
\end{subequations}

Therefore, the KKT conditions of both problems match upon convergence of the modifiers (represented by $l=\infty$) if and only if:

\begin{subequations}
\begin{align}
\begin{split} \label{subeq:kkt:grad}
& \frac{\partial \L_p}{\partial \theta}(\theta^*) = \frac{\partial \L_m}{\partial \theta}(\theta^*) = {\bf{0}}
\end{split} \\
\begin{split} \label{subeq:kkt:base}
& F^\theta_p(\theta^*) = F_m^\theta(\theta^*)+(\Lambda_\infty^\theta)^T \theta^* + \epsilon_\infty.
\end{split}
\end{align}
\end{subequations}

To meet \eqref{subeq:kkt:grad}, we need to set the first order modifiers $\Lambda^\theta$ so that:

\begin{equation*}
 \frac{\partial F^\theta_p}{\partial \theta}(\theta^*) = \frac{\partial F^\theta_m}{\partial \theta}(\theta^*) + \Lambda_\infty^\theta
\end{equation*}

Thus, the optimal modifiers $(\lambda^x_\infty, \lambda^u_\infty)$ must be computed as:

\begin{equation} \label{eq:modifiers:lambda:opt}
\Lambda_\infty^\theta =  \begin{bmatrix} \lambda_\infty^x & \lambda_\infty^u	\end{bmatrix}^T = \frac{\partial F_p^\theta}{\partial \theta}(\theta^*) - \frac{\partial F_m^\theta}{\partial \theta}(\theta^*).
\end{equation}

To converge to the optimal modifiers, we follow an update policy similar to the one proposed in \cite{marchetti2009modifier}. Let $\theta^{\text{drto}}_l$ be the solution of the DRTO \eqref{eq:drto:ma} at iteration $l$, then the modifiers at iteration $l+1$ are calculated as:

\begin{equation} \label{eq:modifiers:lambda}
\Lambda_{l+1}^\theta =  \begin{bmatrix} \lambda_{l+1}^x & \lambda_{l+1}^u	\end{bmatrix}^T = \frac{\partial F_p^\theta}{\partial \theta}(\theta^{\text{drto}}_l) - \frac{\partial F_m^\theta}{\partial \theta}(\theta^{\text{drto}}_l).
\end{equation}

While the gradients of the model can usually be easily computed with user-defined precision, e.g. by numeric or analytical differentiation (Section \ref{sec:grad:linear}), the gradients of the real system often are cumbersome and rely on noisy measures to carry out estimations. The estimation of such gradients is out of the scope of this paper and the reader is referred to other works such as \cite{franccois2014use,vaccari2021offset}.

Given the modifiers $\Lambda_\infty^\theta$ computed in \eqref{eq:modifiers:lambda:opt}, in order to meet \eqref{subeq:kkt:base}, the modifier $\epsilon_\infty$ must be set as:

\begin{equation} \label{eq:modifiers:eps:opt}
\epsilon_\infty =  F^\theta_p(\theta^*) - \left( F_m^\theta(\theta^*)+(\Lambda_\infty^\theta)^T \theta^* \right).
\end{equation}

Applying an update like the one from \eqref{eq:modifiers:lambda}, we get to the following update for $\epsilon_l$:

\begin{equation} \label{eq:modifiers:eps}
\epsilon_{l+1} =  F^\theta_p(\theta^{\text{drto}}_l) - \left( F_m^\theta(\theta^{\text{drto}}_l)+(\Lambda_l^\theta)^T \theta^{\text{drto}}_l \right).
\end{equation}

\subsection{Gradients of a linear model} \label{sec:grad:linear}

In this section, we derive the analytical expression for the gradients of a linear model.

Given the discrete-time linear model

\begin{equation} \label{eq:lineal:model}
\hat{x}_{k+1}= f_{m,k}(x_k,u_k) =A_k \hat{x}_k + B_k u_k,
\end{equation}
we have that

\begin{equation}
\hat{{\bf{x}}}_T = F_m(x_0, {\bf{u}}_T) = \mathcal{F}^x x_0 + \mathcal{F}^u {\bf{u}}_T,
\end{equation}
or equivalently

\begin{equation}
\hat{{\bf{x}}}_T = F_m^\theta(\theta) = \begin{bmatrix} \mathcal{F}^x & \mathcal{F}^u \end{bmatrix} \theta,
\end{equation}
where

\begin{equation}
\begin{aligned}
 \mathcal{F}^x=&\begin{bmatrix} A_0 \\ A_0 A_1 \\ \vdots \\ \prod_{i=0}^{T-1} A_i \end{bmatrix},  \quad
 \mathcal{F}^u=&\begin{bmatrix} B_0 & & & \\ A_1 B_0 & B_1 & &  \\ \vdots & \vdots  &\ddots & \\ (\prod_{i=1}^{T-1 } A_i) B_0 & (\prod_{i=2}^{T-1} A_i) B_1 & \ldots & B_{T-1} \end{bmatrix}.
\end{aligned}
\end{equation}

Therefore, the gradients of the model are constant and can be explicitly computed as $\frac{\partial F_m^\theta}{\partial \theta}= \begin{bmatrix} \mathcal{F}^x & \mathcal{F}^u \end{bmatrix}^T$.

\section{MPC for periodic operation} \label{sec:mpc}

Model predictive controllers (MPCs) are one of the multiple choices for the bottom layer of the two-layer scheme introduced in Section \ref{sec:problem:formulation}, often refered to as advanced control. Its objective is to calculate the control sequence that takes the system from its current state $z_j$ to the reference given by the STTO.

Contrary to the DRTO, the MPC generally has a more local and fast nature, which may make its correspondent real system different from the one presented in \eqref{eq:real:sys}. Let the local system be defined as

\begin{equation} \label{eq:real:sys:mpc}
z_{j+1} = f^\text{mpc}_{p,j}(z_j,v_j),
\end{equation}
where $z_j \in \R^{n_z}$ and $v_j \in \R^{n_v}$ represent respectively the states and inputs of the local system at time $j$, and $f^\text{mpc}_{p,j} : \R^{n_z \times n_v} \to \R^{n_z}$ represents the dynamics of the real system at time $j$. Note that the local system is parameterized by $j$ to indicate that the discretization time of the local system ($t_N$ seconds) is generally different than that of the global system parameterized by $k$ ($t_T$ seconds). Therefore, the system is periodic with period $t_N L = t_T T$ seconds.

The  MPC solves at each time step $j$ an optimization problem to calculate the optimal sequence of control inputs. Given a reference at time $j$ $({\bf{\hat{z}}}_{N,j}^{\text{ref}}, {\bf{v}}_{N,j}^{\text{ref}})$, we use an offset free MPC formulation based on \cite{muske2002disturbance,pannocchia2003disturbance}.

Let the model of the MPC local system \eqref{eq:real:sys:mpc} be defined as

\begin{equation}
\hat{z}_{j+1} =  f^\text{mpc}_m(\hat{z}_j,v_j,d_j),
\end{equation}
where $f^\text{mpc}_m$ is usually a linear system to allow fast MPC implementations and $d_j \in \R^{n_z}$ is the so-called disturbance at time $j$. In contrast to the local system \eqref{eq:real:sys:mpc}, the model $f^\text{mpc}_m$ is time-invariant and its dependence of time comes through the disturbances $d_j$. Local constraints are also considered, but for the sake of simplicity, only as box constraints on the inputs $v_j$. More general constraints require robust formulations of the MPC to guarantee recursive feasibility \cite{limon2010robust,limon2008design,bemporad2007robust} and are out of the scope of this work.

The sequence of disturbances $d_j$ is periodic over the periodic horizon $L$ and are updated in such a way that, upon convergence, 

\begin{equation} \label{eq:d:periodic}
f^\text{mpc}_m(z_j,v_j,d_j) = f^\text{mpc}_{p,j}(z_j,v_j).
\end{equation}

Now we present a simple way to estimate the disturbances

\begin{equation} \label{eq:disturbance:future}
d_{j+L} = d_j + K^d \left( f^\text{mpc}_{p,j}(z_j,v_j)  - f^\text{mpc}_m(z_j,v_j,d_j) \right),
\end{equation}
where $K^d \in \R^{n_z \times n_z}$ is a filtering matrix and the matrix $K_d$ is stable.

Given the current state $z_j$ and the sequence of future disturbances ${\bf{d}}_N = \begin{bmatrix} d_j &d_{j+1} &\hdots &d_{j+N-1}   \end{bmatrix}$, the offset-free periodic MPC at time step $j$ is formulated as follows:

\begin{equation} \label{eq:mpc}
\begin{aligned}
({\bf{z}}^*,{\bf{v}}^*) & = \\
\argmin_{{\bf{\hat{z}}}_N,{\bf{v}}_N} \quad & \ell^{\text{mpc}}({\bf{\hat{z}}}_N,{\bf{v}}_N,{\bf{\hat{z}}}_{N,j}^{\text{ref}}, {\bf{v}}_{N,j}^{\text{ref}}) \\
\text{s.t.} \quad & \hat{z}_{i+1} = f^{\text{mpc}}_m(\hat{z}_i,v_i,d_{j+i}), \; \text{for all } i=0, 1, \ldots, N-1 \\
& v^L \leq v_i \leq v^U,  \; \text{for all } i=0, 1, \ldots, N-1 \\
& \hat{z}_0 = z_j,
\end{aligned}
\end{equation}
where $\ell^{\text{mpc}}$ is a cost function that penalizes the distance between the reference sequences of states and inputs, and ${\bf{\hat{z}}}_N = \begin{bmatrix} \hat{z}_1 &\hat{z}_2 & \hdots &\hat{z}_N  \end{bmatrix}, {\bf{v}}_N = \begin{bmatrix} v_0 &v_1 & \hdots &v_{N-1} \end{bmatrix}$. The current local state $z_j$ is considered known. The MPC follows a receding horizon scheme, which means that only the first computed input $v^*_0$ is applied to the system at each iteration of the MPC.

\section{Steady trajectory target optimization (STTO)} \label{sec:stto}

The solution of the DRTO presented in Section \ref{sec:pma} leads to the reference trajectory $r_e^{\text{drto}}$. The objective of the STTO is to transform this reference trajectory into a valid target $( {\bf{\hat{z}}}^{\text{ref}}_{L,j}, {\bf{v}}^{\text{ref}}_{L,j} )$ for the MPC defined in Section \ref{sec:mpc}, i.e. one feasible for the MPC constraints.

The first step is to match the time scale of the DRTO ($t_T$) with that of the MPC ($t_N$). Usually, the DRTO works with longer time steps than the MPC ($t_T > t_N$). Therefore, one must transform the reference given by the DRTO into one with the same time scale of the MPC. The reference given by the DRTO spans a total duration of $T t_T$ seconds. To transform it into the time scale of the MPC, just divide it into segments of $t_N$ seconds and check which value of the reference trajectory ${\bf{r}}_e^{\text{drto}}$ corresponds to each segment. To avoid dealing with segments that comprise two or more values, we assume that the DRTO sampling time $t_T$ is a multiple of the MPC sampling time $t_N$. The new reference with time scale $t_N$ is denoted ${\bf{r}}_e^{\text{stto}}$ and its length is $L$. Then, this reference is shifted to match the current time step $j$.

Let $\ell^{\text{stto}} : \R^{L n_z \times L n_v \times L n_r} \to \R$ be a function that penalizes the distance between the trajectories of states and inputs of the MPC (${\bf{\hat{z}}}, {\bf{v}}$) and the reference trajectory (${\bf{r}}_e^{\text{stto}})$. Then, the STTO problem at step $j$ can be formulated as:

\begin{equation}
\begin{aligned}
( {\bf{\hat{z}}}^{\text{ref}}_{L,j}, {\bf{v}}^{\text{ref}}_{L,j} ) =& \\
\argmin_{{\bf{\hat{z}}},{\bf{v}}} \quad & \ell^{\text{stto}} ({\bf{\hat{z}}},{\bf{v}},{\bf{r}}_e^{\text{stto}}) \\
\text{s.t.} \quad & \hat{z}_{i+1} = f^{\text{mpc}}_m(\hat{z}_i,v_i,d_{j+i}), \; \text{for all } i=0, 1, \ldots, L-1 \\
& v^L \leq v_i \leq v^U, \; \text{for all } i=0,1, \ldots, L-1 \\
& \hat{z}_L = \hat{z}_0,
\end{aligned}
\end{equation}
where ${\bf{d}}_L = \begin{bmatrix} d_j &d_{j+1} &\hdots &d_{j+L-1}   \end{bmatrix}$ are the disturbances estimated in \eqref{eq:disturbance:future}.

Finally, given that $t_T > t_N$, then we have that $L>N$ and therefore the obtained reference $( {\bf{\hat{z}}}^{\text{ref}}_{L,j}, {\bf{v}}^{\text{ref}}_{L,j} )$ must be trimmed to match the control horizon $N$. This new reference trajectory will be referred to as $({\bf{\hat{z}}}_{N,j}^{\text{ref}}, {\bf{v}}_{N,j}^{\text{ref}})$ and constitutes a valid reference for the MPC layer which is guaranteed to be feasible for the MPC constraints.

As commented in Section \ref{sec:problem:formulation}, the STTO layer needs to be computed before every MPC iteration to guarantee that the reference trajectory is feasible for the most recent values of disturbances. Since the DRTO reference ${\bf{r}}_e^{\text{drto}}$ is periodic, endless shifting is possible and the STTO/MPC loop can always control the system to the newest reference trajectory.

In the next section, we detail the full algorithms for the DRTO/STTO and the MPC layers.

\section{Periodic Modifier-Adaptation Algorithm} \label{sec:algorithm}

In this section, we go through the full periodic modifier-adaptation algorithm. As commented in Section \ref{sec:problem:formulation}, the full scheme can be splitted into two main layers with different time scales. Algorithm \ref{alg:drto} goes through the DRTO layer, which uses a time scale of $t_D$ seconds. It shows how to compute the reference for the STTO and MPC layers. Besides, Algorithm \ref{alg:mpc} details the STTO and MPC layers which work with a shorter time scale $t_N$ and calculates every control signal that is applied to the system. Both algorithms run in parallel as long as automatic control of the system is required. Upon convergence to the optimal predicted trajectory, the controlled system is guaranteed to reach an optimal periodic behaviour.

\begin{algorithm}
Init
\begin{enumerate}[(i)]
\item Initialize $l=0$ and the modifiers $\hat \lambda^x_0, \hat \lambda^u_0, \hat \epsilon_0$ to zero.
\end{enumerate}

Loop
\begin{enumerate}[(i)]
\setcounter{enumi}{1}
\item Given the modifiers $\hat \lambda^x_l, \hat \lambda^u_l, \hat \epsilon_l$, compute the optimal trajectory with the DRTO defined in \eqref{eq:drto:ma} and \eqref{eq:drto:re} and obtain $ {\bf{r}}_e^{\text{drto}} $.
\item \label{alg:drto:mpc} Pass ${\bf{r}}_e^{\text{drto}}$ as a reference to the STTO and MPC Algorithm.
\item \label{alg:drto:ma} Estimate the gradients of the model and the real system and update the modifiers according to \eqref{eq:modifiers:lambda} and \eqref{eq:modifiers:eps}:
\begin{equation*}
\begin{aligned}
 \lambda_{l+1}^x= & \frac{\partial F_p}{\partial x_0}\Big\rvert_{(x^{\text{drto}}_0,{\bf{u}}_{T}^{\text{drto}})}-\frac{\partial F_m}{\partial x_0}\Big\rvert_{(x^{\text{drto}}_0,{\bf{u}}_{T}^{\text{drto}})} \\
 \lambda_{l+1}^u= & \frac{\partial F_p}{\partial {\bf{u}}_T}\Big\rvert_{(x^{\text{drto}}_0,{\bf{u}}_{T}^{\text{drto}})}-\frac{\partial F_m}{\partial {\bf{u}}_T}\Big\rvert_{(x^{\text{drto}}_0,{\bf{u}}_{T}^{\text{drto}})} \\
\epsilon_{l+1} = & F_p(x^{\text{drto}}_0,{\bf{u}}_{T}^{\text{drto}}) - \\
& \left(F_m(x^{\text{drto}}_0,{\bf{u}}_{T}^{\text{drto}}) + \hat{\lambda}_l^x x_0^{\text{drto}} + \hat{\lambda}_l^u {\bf{u}}_T^{\text{drto}} \right)
\end{aligned}
\end{equation*}

\item Wait until next iteration and update $l=l+1$.
\end{enumerate}
End Loop
\caption{DRTO algorithm (executed each $t_D$ seconds)}
 \label{alg:drto}
\end{algorithm}

\begin{algorithm}
Init
\begin{enumerate}[(i)]
\item Initialize disturbances of the first period to zero, i.e. $d_j=0$, where $j=0,1, \ldots, L-1$, and set $j=0$.
\end{enumerate}
Loop
\begin{enumerate}[(i)]
\setcounter{enumi}{1}
\item Get the current state $z_j$.
\item \label{alg:mpc:stto} Given $ {\bf{r}}_e^{\text{drto}} $ from the DRTO Algorithm, shift it to match the current time step $j$ and solve the STTO layer detailed in Section \ref{sec:stto} to obtain $( {\bf{\hat{z}}}^{\text{ref}}_{N,j}, {\bf{v}}^{\text{ref}}_{N,j} )$.
\item \label{alg:mpc:mpc} Given the reference trajectory $( {\bf{\hat{z}}}^{\text{ref}}_{N,j}, {\bf{v}}^{\text{ref}}_{N,j} )$, compute the MPC control input $v_j$ from \eqref{eq:mpc}.
\item Apply input $v_j$ to the local system.
\item Estimate the disturbance for the next period $d_{j+L}$ using \eqref{eq:disturbance:future}.
\item Wait until next iteration and update $j=j+1$.
\end{enumerate}
End Loop
\caption{STTO and MPC algorithm (executed each $t_N$ seconds)}
 \label{alg:mpc}
\end{algorithm}

\begin{remark} \label{rem:filter}
An optional filtering of the modifiers can be performed after step \eqref{alg:drto:ma} of Algorithm \ref{alg:drto}, filtering the modifiers can influence the speed and stability of the convergence process:
\begin{equation*}
\begin{aligned}
\tilde{\lambda}^x_{l+1}&= K^x \lambda_{l+1}^x	+ (\Id-K^x) \tilde{\lambda}_{l}^x \\
\tilde{\lambda}^u_{l+1}&= K^u \lambda_{l+1}^u	 + (\Id-K^u) \tilde{\lambda}_{l}^u \\
\epsilon_{l+1}&= K^\epsilon \epsilon_{l+1} + (\Id-K^\epsilon) \epsilon_{l},
\end{aligned}
\end{equation*}
where $K^x, K^u$ and $K^\epsilon$ are stable filter matrices and $( \; \tilde{ } \; )$ notation refers to the filtered modifiers.
\end{remark}

\begin{remark} \label{rem:layer:independence}
The formulations of the STTO and the MPC in steps \eqref{alg:mpc:stto} and \eqref{alg:mpc:mpc} of Algorithm \ref{alg:mpc} are not fixed and other options besides the ones presented on this paper are equally acceptable.
\end{remark}

In the next section, we show the performance of P-MA in the quadruple tank benchmark.

\section{Illustrative example} \label{sec:example}

In this section we show the performance of the periodic modifier-adaptation formulation of the DRTO introduced in this paper. For the sake of clarity, we omit the STTO and MPC layers and show how the reference trajectory computed by the P-MA formulation of the DRTO converges to the optimal periodic trajectory.

To study the performance of the proposed approach, we test it against a periodic version of the quadruple tank process. This benchmark first proposed in \cite{johansson2000quadruple} has been widely used to test different controllers \cite{alvarado2011comparative}. The quadruple-tank system scheme is shown in Figure \ref{fig:4tanks:scheme} and consists of four interconnected tanks that share water according to the following physical equations:

\begin{figure}
\centering
\includegraphics[width=0.45\linewidth]{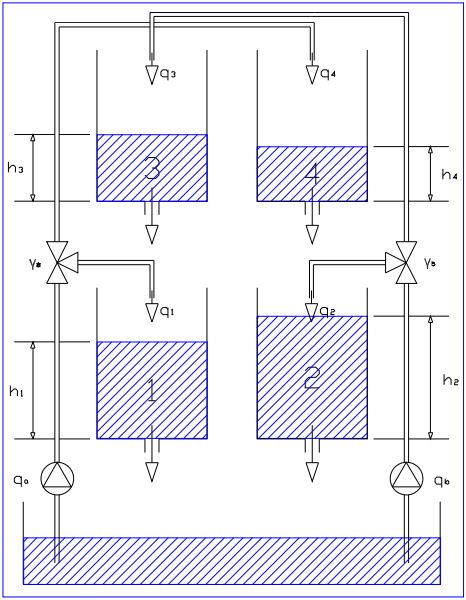}
\caption{Quadruple-tank system diagram, reproduced from \cite{alvarado2011comparative}.}
\label{fig:4tanks:scheme}
\end{figure}

\begin{equation}
\begin{aligned}
& S \frac{dh_{1}}{dt}=-a_{1}\sqrt{2gh_{1}}+a_{3}\sqrt{2gh_{3}}+\frac{\gamma_{a}q_{a}}{3600}\\
& S \frac{dh_{2}}{dt}=-a_{2}\sqrt{2gh_{2}}+a_{4}\sqrt{2gh_{4}}+\frac{\gamma_{b}q_{b}}{3600}\\
& S \frac{dh_{3}}{dt}=-a_{3}\sqrt{2gh_{3}}+(1-\gamma_{b})\frac{q_{b}}{3600}\\
& S \frac{dh_{4}}{dt}=-a_{4}\sqrt{2gh_{4}}+(1-\gamma_{a})\frac{q_{a}}{3600}.
\end{aligned}
\label{eq:4tanks}
\end{equation}

And are subject to the following box constraints:

\begin{equation} \label{eq:4tanks:cons}
{\bf{h}}_\text{min} \leq \begin{bmatrix} h_1 \\ h_2 \\ h_3 \\ h_4 \end{bmatrix} \leq {\bf{h}}_\text{max}		\qquad		{\bf{q}}_\text{min} \leq \begin{bmatrix} q_a \\ q_b \end{bmatrix} \leq {\bf{q}}_\text{max}.
\end{equation}

The quadruple tank process has some relevant properties:
\begin{itemize}
\item It presents large coupling between its subsystems.
\item It dynamics are nonlinear.
\item States can be measured.
\item States and inputs are hard constrained.
\item Its real gradients can be analytically computed with the physical equations.
\end{itemize}

We use a compact notation to define the parameters of the plant:

$$\bf{a}=\begin{bmatrix} a_1 \\ a_2 \\ a_3 \\ a_4 \end{bmatrix}, \bf{x} = \begin{bmatrix} h_1 \\ h_2 \\ h_3 \\ h_4 \end{bmatrix}, \bf{u} = \begin{bmatrix} q_a \\ q_b \end{bmatrix}, \boldsymbol{\gamma} = \begin{bmatrix} \gamma_a \\ \gamma_b \end{bmatrix},$$
where water levels ${\bf{x}}$ corresponds to the states and water flows ${\bf{u}}$ to the inputs of the system. Information about the parameters is collected in Table \ref{tab:4tanks:parameters} and \eqref{eq:4tanks:gamma}. The periodic nature of the system is induced through parameter $\boldsymbol{\gamma}$, whose cycle is shown in \eqref{eq:4tanks:gamma}, where each column represent a constant value of $\boldsymbol{\gamma}$ for $t_T= 3600$ seconds. Therefore, the plant is periodic with period $T=7$ hours.

\begin{table}
\caption{Parameters of the plant}
\begin{center}
\begin{tabular}{c| c | c | c} 
  & Value & Unit & Description \\ \hline \hline
$S$ &  0.03 & m$^2$ & Cross-section of the tanks \\ 
$\bf{a}$ & $\begin{bmatrix} 1.31 \\ 1.51 \\ 0.927 \\ 0.882	\end{bmatrix} e^{-4}$ & m$^2$ & Discharge constants\\ 
${\bf{h}}_\text{max}$ & $\begin{bmatrix} 1.36 \\ 1.36 \\ 1.30 \\ 1.30 \end{bmatrix}$ & m & Maximum water level\\ 
${\bf{h}}_\text{min}$ & $\begin{bmatrix} 0.2 \\ 0.2 \\ 0.2 \\ 0.2 \end{bmatrix}$ & m & Minimum water level\\ 
${\bf{q}}_\text{max}$ & $\begin{bmatrix} 3.6 \\ 4.0 \end{bmatrix}$ & m$^3$/h & Maximum water flow\\ 
${\bf{q}}_\text{min}$ & $\begin{bmatrix} 0 \\ 0 \end{bmatrix}$ & m$^3$/h & Minimum water flow\\ 
$g$ & 9.81 & m/s$^2$ & Gravity acceleration \\
\end{tabular}
\end{center}
\label{tab:4tanks:parameters}
\end{table}

\begin{equation} \label{eq:4tanks:gamma}
\boldsymbol{\gamma}_{\text{cycle}}=\begin{bmatrix}0.3& 0.4 &  0.5& 0.7 &0.6  &0.4  &  0.2\\
    0.6 &   0.5 &   0.4&    0.2 &  0.3 &   0.5 &   0.7\end{bmatrix}.
\end{equation}

The model of the system \eqref{eq:4tanks} is a discrete linear model with discretization time set to $5$ seconds and linearized at the point:

\begin{equation*}
\begin{aligned}
x_0=\begin{bmatrix}0.7293 \\ 0.8102 \\ 0.6594 \\ 0.9408 \end{bmatrix}, \quad u_0 = \begin{bmatrix}1.948 \\ 2.00 \end{bmatrix}, \quad \boldsymbol{\gamma}_0=\begin{bmatrix} 0.3 \\ 0.4 \end{bmatrix}
\end{aligned}
\end{equation*}

Therefore, the model can be written as:

\begin{equation}
\begin{aligned}
x_{k+1} = & \begin{bmatrix} 0.945 &   0  &  0.040 &    0\\
         0  &  0.940&   0 &  0.032\\
         0  & 0   & 0.959  &  0\\
         0  & 0   & 0 & 0.967\end{bmatrix} (x_k - x_0) + \\
         &  \begin{bmatrix}0.0135&  0.0006\\
    0.0005  & 0.0180\\
         0 & 0.0272\\
    0.0319 & 0\end{bmatrix} (u_k - u_0) + x_0.
\end{aligned}
\end{equation}

At every time step, the system is subject to the box constraints on inputs and states from \eqref{eq:4tanks:cons}, i.e.

\begin{equation} \label{eq:4tanks:cons:model}
{\bf{h}}_\text{min} \leq {\bf{x}} \leq {\bf{h}}_\text{max}		\qquad		{\bf{q}}_\text{min} \leq {\bf{u}} \leq {\bf{q}}_\text{max}.
\end{equation}

Given the economic parameters $c=1$ and $p=20$, the economic cost of operating the plant at each discrete time step is given by

$$ \phi(x_k,u_k) = (q_a^2 + c q_b^2) + p \frac{0.012}{S (h_1 + h_2)}.$$

We apply Algorithm \ref{alg:drto} to compute the optimal periodic trajectory for the system. The control process, i.e. steps \eqref{alg:mpc:stto} and \eqref{alg:drto:mpc} and Algorithm \ref{alg:mpc}, is omited for the sake of clarity. The periodic constant is taken as $T=7$ and the optional filtering of the modifiers proposed in Remark \ref{rem:filter} has not been taken into account.

The integration of the real process as well as the computation of the optimal trajectory and the P-MA DRTO reference trajectory have been computed using the CasADi optimization tool in Matlab \cite{Andersson2018}. The gradients of the real process have been computed using numerical differentiation on the real system \eqref{eq:4tanks}, while those of the linear model have been computed using the results from Section \ref{sec:grad:linear}.

\begin{figure}
\centering
\includegraphics[width=1\linewidth]{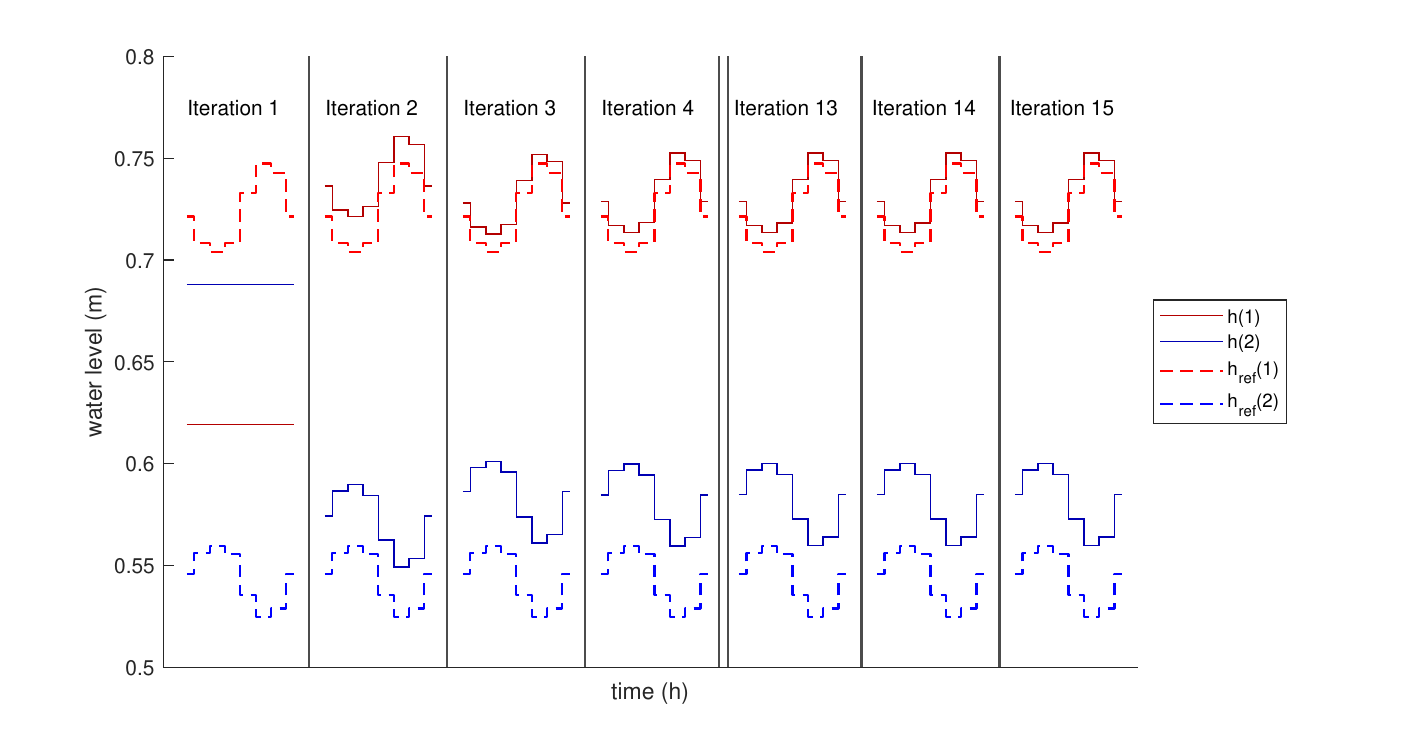}
\caption{Sequence of states computed by the DRTO with the zeroth order modifier in different iterations (solid lines) vs optimal sequence of states (dashed lines).}
\label{fig:4tanks:states:noma}
\end{figure}

\begin{figure}
\centering
\includegraphics[width=1\linewidth]{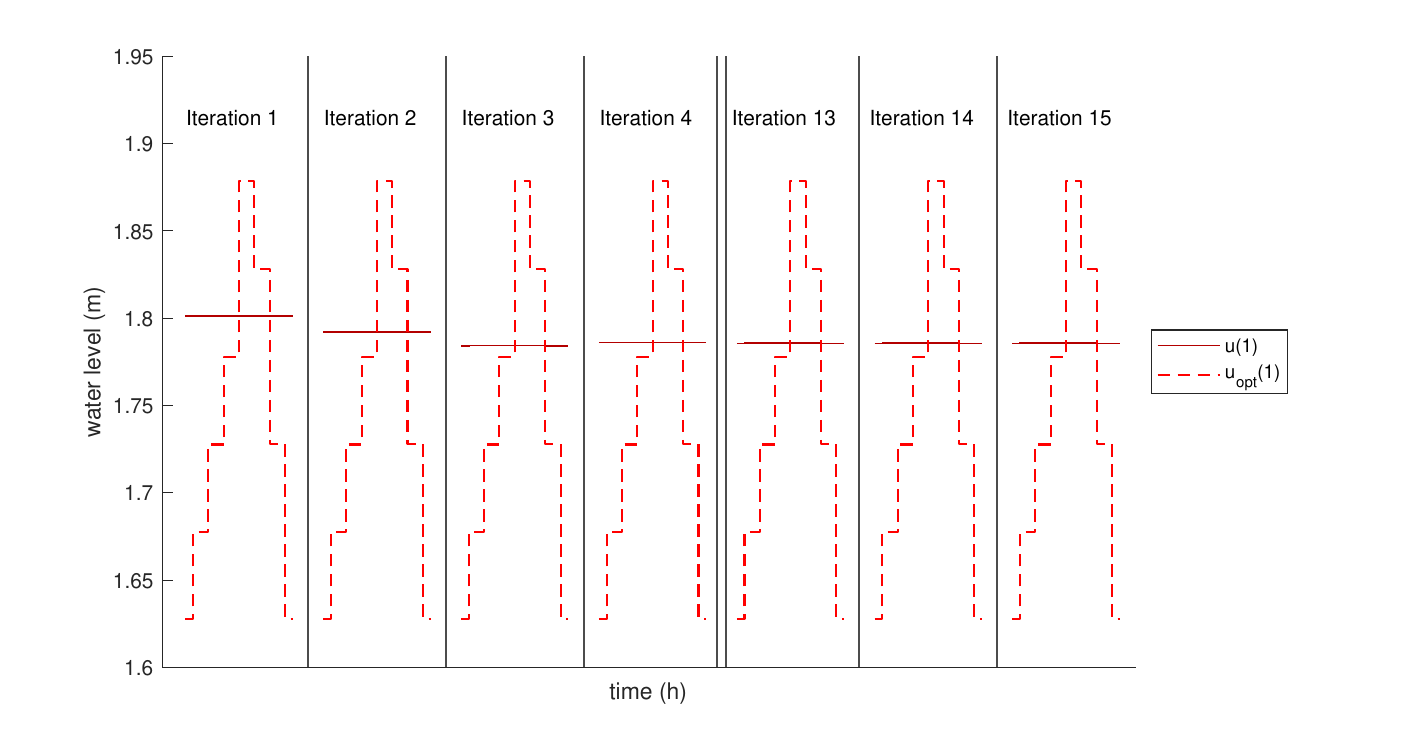}
\includegraphics[width=1\linewidth]{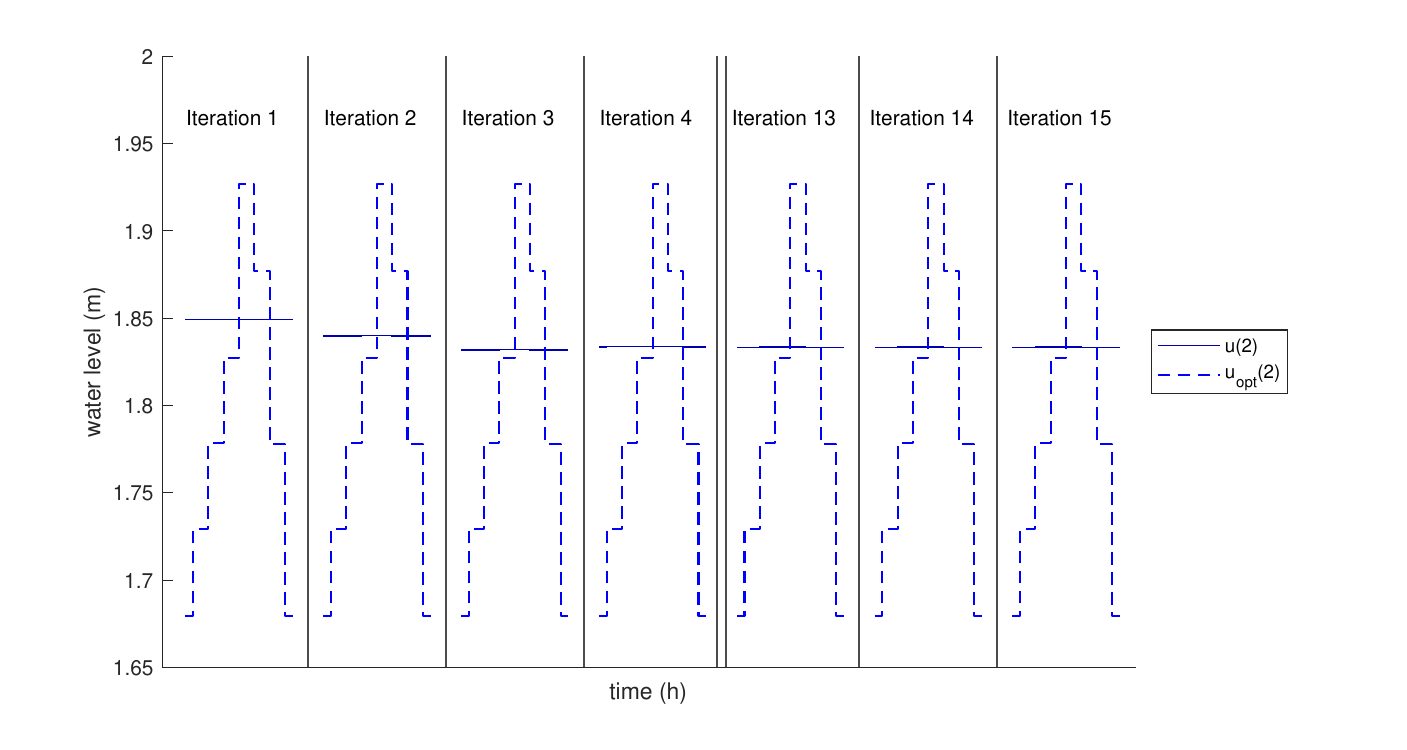}
\caption{Sequence of inputs calculated by the DRTO with the zeroth order modifier in different iterations (solid lines) vs optimal sequence of inputs (dashed lines).}
\label{fig:4tanks:inputs:noma}
\end{figure}

Figures \ref{fig:4tanks:states:noma} and \ref{fig:4tanks:inputs:noma} show the optimal trajectory computed by the DRTO with no first order modifiers. The inclusion and convergence of the zeroth order modifier $\epsilon_l$ guarantees that the predicted trajectory matches the response of the real system. However, the lack of first order modifiers entails that, upon convergence, the computed input sequence is not optimal. Moreover, since the KKT conditions of this DRTO does not change with time, the sequence of inputs predicted at each iteration is constant over time.

Notice how in Iteration 1 (Figures \ref{fig:4tanks:states:noma} and \ref{fig:4tanks:inputs:noma}), all the modifiers are set to zero and the optimal predicted behaviour is a single steady state. This is due to the time invariant model used by the DRTO, which differs vastly from the real periodic behaviour of the system.

\begin{figure}
\centering
\includegraphics[width=1\linewidth]{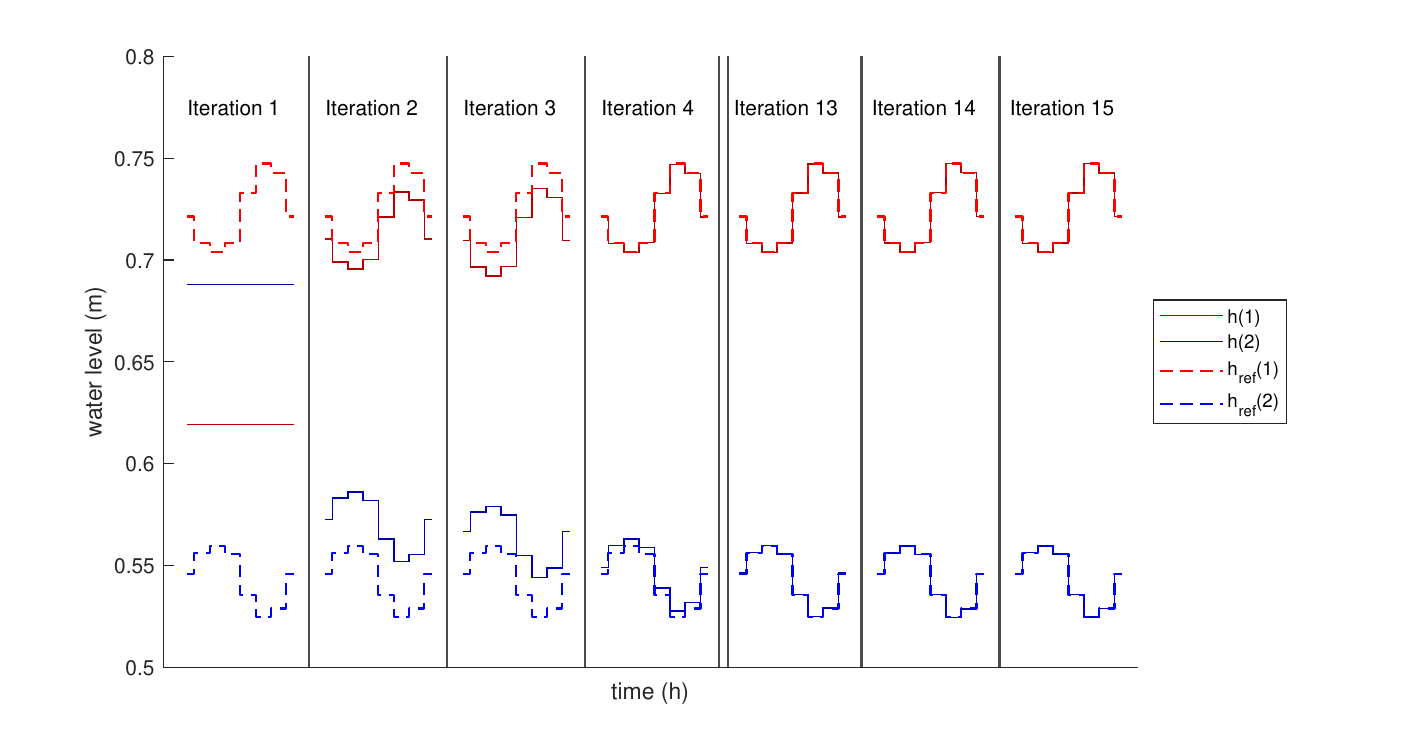}
\caption{Sequence of states computed by the P-MA DRTO in different iterations (solid lines) vs optimal sequence of states (dashed lines).}
\label{fig:4tanks:states:ma}
\end{figure}

\begin{figure}
\centering
\includegraphics[width=1\linewidth]{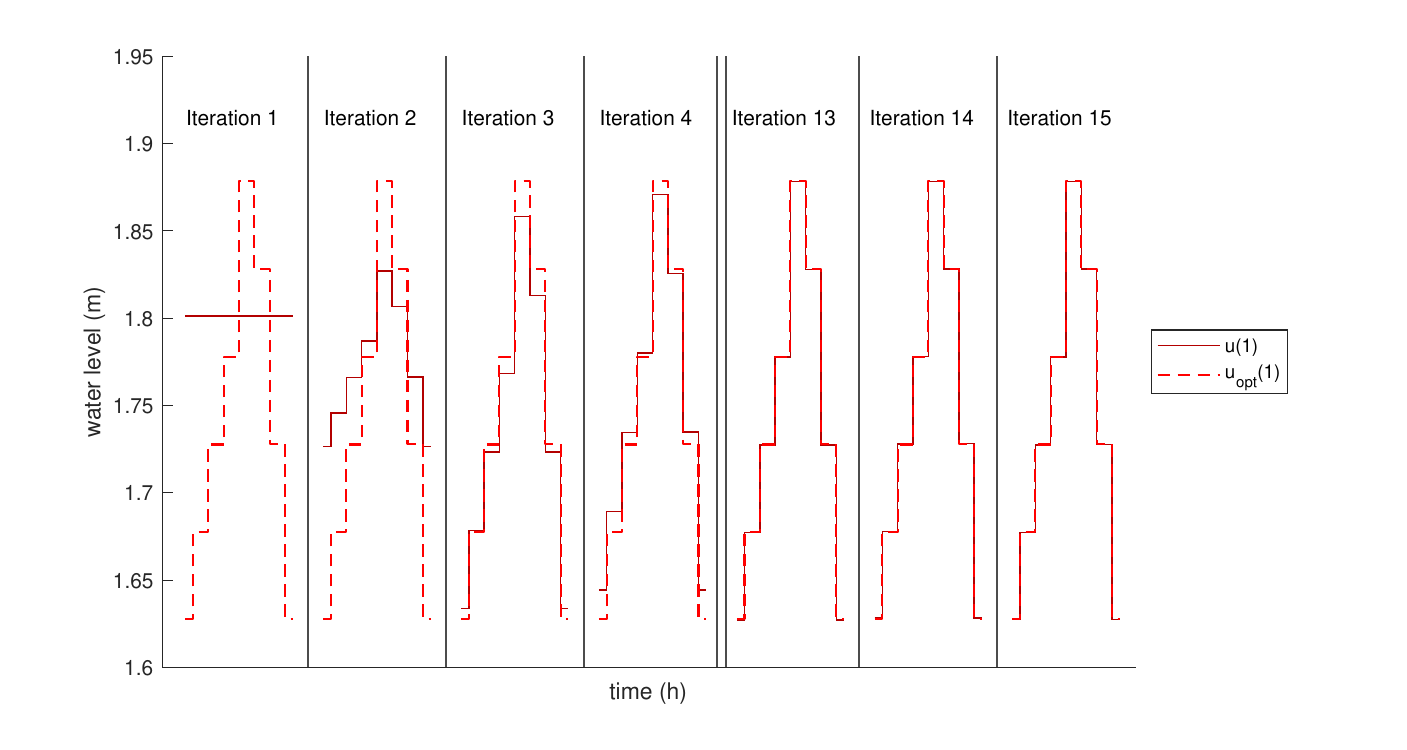}
\includegraphics[width=1\linewidth]{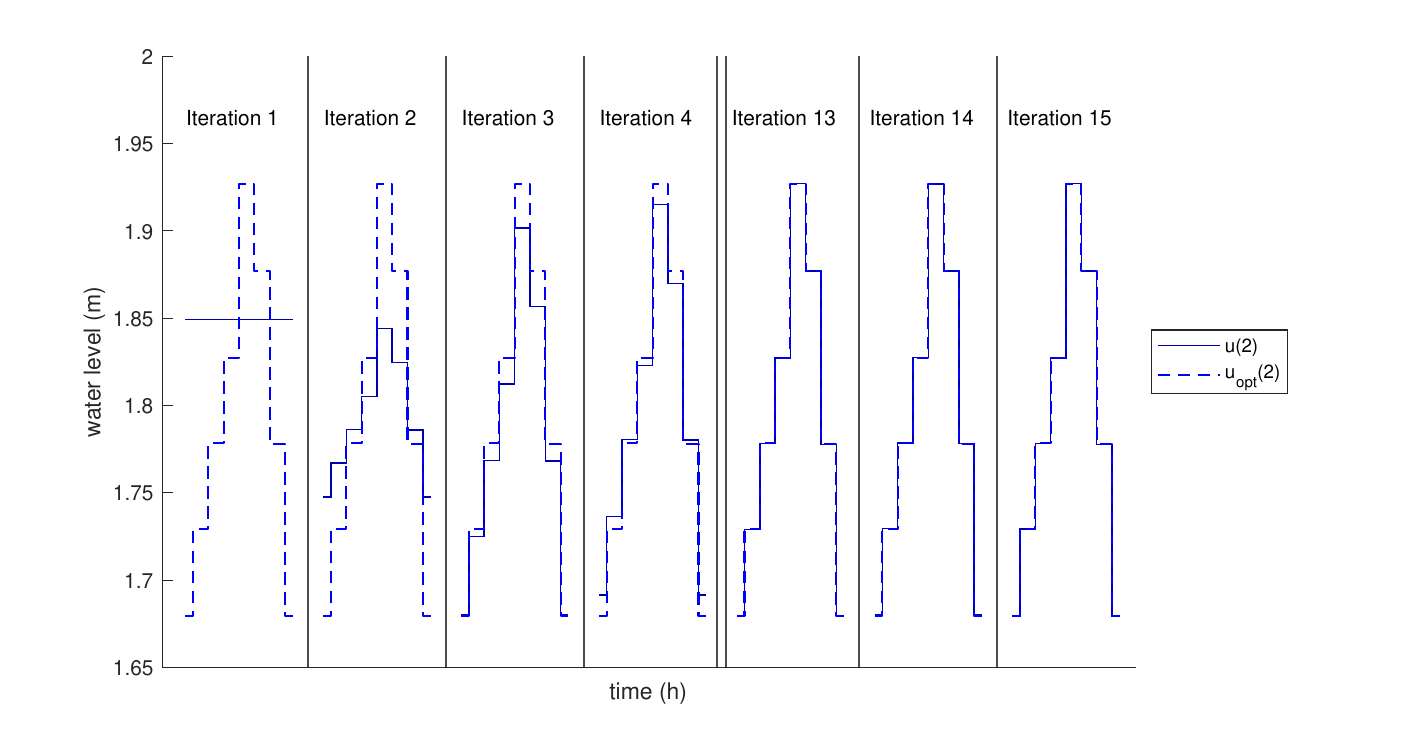}
\caption{Sequence of inputs calculated by the P-MA DRTO in different iterations (solid lines) vs optimal sequence of inputs (dashed lines).}
\label{fig:4tanks:inputs:ma}
\end{figure}

Figure \ref{fig:4tanks:states:ma} shows how the P-MA DRTO achieves convergence to the optimal sequence of states, and Figure \ref{fig:4tanks:inputs:ma} shows that this convergence is also achieved with the optimal sequence of inputs. After 15 iterations, the sequences of states and inputs computed by the P-MA DRTO are sufficiently close to the optimal sequences.

\section{Conclusions} \label{sec:conclusions}

In this paper we have presented a periodic scheme which is able to control a system given inaccurate models of it.

First, the optimal trajectory is computed by the P-MA DRTO, which uses information of the real system to modify the dynamic real time optimization layer with first and zeroth order modifiers so that, upon convergence of these modifiers, the KKT conditions of the P-MA DRTO converge to the optimal operation of the real plan. Then, a steady trajectory target optimization (STTO) translates the trajectory computed by the P-MA DRTO into a feasible reference for the MPC. Finally, the MPC layer uses a disturbance estimator to adapt its model to the real system and converge to the optimal reference.

The proposed P-MA DRTO has been tested against the quadruple tank process, showing that its solution does converge to the optimal periodic behaviour given that the gradients of the real system are computed with enough accuracy.

\bibliographystyle{elsarticle-num}
\bibliography{BibNew}

\end{document}